\documentclass{amsproc}
\usepackage{amssymb}

\newtheorem{theorem}{Theorem}[section]
\newtheorem{lemma}[theorem]{Lemma}
\theoremstyle{definition}
\newtheorem{definition}[theorem]{Definition}
\newtheorem{example}[theorem]{Example}
\newtheorem{proposition}[theorem]{Proposition}
\newtheorem{corollary}[theorem]{Corollary}

\theoremstyle{remark}

\numberwithin{equation}{section}

\begin{document}

\title[Symbolic Rees algebras of edge ideals]{Combinatorics of
symbolic Rees algebras of edge ideals of clutters}

\author[Mart\'{\i}nez-Bernal]{Jos\'e Mart\'{\i}nez-Bernal}
\address{
Departamento de
Matem\'aticas\\
Centro de Investigaci\'on y de Estudios Avanzados del
IPN\\
Apartado Postal
14--740 \\
07000 Mexico City, D.F. } \email{jmb@math.cinvestav.mx}
\thanks{The second author was Partially supported by 
COFAA-IPN and SNI. The third author was partially supported by CONACyT 
grant 49251-F and SNI}

\author[Renter\'\i a-M\'arquez]{Carlos Renter\'\i a-M\'arquez}
\address{
Departamento de Matem\'aticas\\
Escuela Superior de F\'\i sica y
Matem\'aticas\\
Instituto Polit\'ecnico Nacional\\
07300 Mexico City, D.F.}
\email{renteri@esfm.ipn.mx}

\author[Villarreal]{Rafael H. Villarreal}
\address{
Departamento de
Matem\'aticas\\
Centro de Investigaci\'on y de Estudios
Avanzados del
IPN\\
Apartado Postal
14--740 \\
07000 Mexico City, D.F.
}
\email{vila@math.cinvestav.mx}

\dedicatory{This paper is dedicated to Wolmer Vasconcelos.}

\begin{abstract} Let $\mathcal{C}$ be a clutter and let $I$ be its
 edge ideal. We present a combinatorial description of the minimal generators of the 
symbolic Rees algebra $R_s(I)$ of $I$. 
It is shown that the minimal generators of $R_s(I)$ are in one to one
correspondence with the indecomposable 
parallelizations of $\mathcal{C}$. From our description some major
results on symbolic Rees algebras of perfect graphs and clutters 
will follow. As a byproduct, we give a method, 
using Hilbert bases, to compute all indecomposable parallelizations of 
$\mathcal{C}$ and all the corresponding vertex covering
numbers. 
\end{abstract}

\maketitle

\section{Introduction}

Let $\mathcal{C}$ be a {\it clutter\/} with finite vertex set 
$X=\{x_1,\ldots,x_n\}$, i.e., $\mathcal{C}$ is a family of 
subsets of $X$, called edges, none 
of which is included in another. The set of
vertices and edges of $\mathcal{C}$ are denoted by $V(\mathcal{C})$
and $E(\mathcal{C})$ respectively. A basic example  
of a clutter is a graph. Let $R=K[x_1,\ldots,x_n]$ be a polynomial ring 
over a field $K$. The {\it edge ideal\/} of $\mathcal{C}$, 
denoted by $I=I(\mathcal{C})$, is the ideal of $R$
generated by all square-free monomials $x_e=\prod_{x_i\in e}x_i$ such 
that $e\in E(\mathcal{C})$. The assignment $\mathcal{C}\mapsto
I(\mathcal{C})$ gives a natural one to one
correspondence between the family of clutters and the family of
square-free monomial ideals.

The {\it blowup algebra\/} studied here is the {\it symbolic Rees
algebra\/} of $I$:
$$
R_s(I)=R\oplus I^{(1)}t\oplus\cdots\oplus I^{(i)}t^i\oplus\cdots\subset R[t],
$$
where $t$ is a new variable and $I^{(i)}$ is the 
$i${\it th\/} symbolic power of $I$. Recall that the 
$i${\it th} {\it symbolic power\/} of
$I$ is defined as
\[
I^{(i)}=S^{-1}I^i\cap R,
\]
where $S=R\setminus\cup_{k=1}^s{\mathfrak p}_i$,
the ideals $\mathfrak{p}_1,\ldots,\mathfrak{p}_s$ are the minimal primes of $I$ 
and $S^{-1}I^i$ is the localization of $I^i$ at $S$. In our situation
the $i${\it th} symbolic power of $I$ can be expressed using systems of
linear inequalities (see Eq.~(\ref{march5-09}) in
Section~\ref{blowupcovers}). Closely related to $R_s(I)$ is---another
blowup algebra---the {\it Rees algebra\/} of $I$:
$$
R[It]=R\oplus It\oplus\cdots\oplus I^{i}t^i\oplus\cdots
=K[\{x_1,\ldots,x_n,x_et\vert\, e\in E(\mathcal{C})\}]\subset R[t].
$$
Blowup algebras are interesting objects of study in algebra and
geometry \cite{Vas}. 

The study of symbolic powers of edge ideals from the point of view of
graph theory and combinatorics was initiated in 
\cite{ITG} and further elaborated on in
\cite{sullivant,perfect}. A breakthrough in this area is the
translation of combinatorial problems (e.g., the Conforti-Cornu\'ejols 
conjecture \cite{cornu-book}, the max-flow min-cut property, the
idealness of a clutter, or the integer rounding property) 
into algebraic problems of blowup algebras of edge ideals
\cite{ainv,roundp,poset,mfmc,reesclu,clutters,chordal}.  

By a result of
Lyubeznik \cite{Lyu3}, $R_s(I)$ is a $K$-algebra of finite type 
generated by a unique minimal finite set of monomials. The
main theorem of this paper is a description---in combinatorial optimization
terms---of this minimal set of generators of $R_s(I)$ 
as a $K$-algebra. Before stating the 
theorem, we need to recall some more terminology and notations.

A subset $C$ of $X$ is called a {\it vertex cover\/} of 
$\mathcal C$ if every edge of $\mathcal C$ contains at least one
vertex of $C$. A vertex cover $C$ is called a {\it minimal 
vertex cover\/} if no proper subset of $C$  is a vertex cover. 
The number of vertices in a minimum vertex 
cover of $\mathcal{C}$, denoted by 
$\alpha_0(\mathcal{C})$, is called the {\it vertex
covering number\/} of $\mathcal{C}$. The dual concept of a vertex
cover is a {\it stable set\/}, i.e., a 
subset $C$ of $X$ is a 
vertex cover of
$\mathcal{C}$  
if and only if $X\setminus C$ is a stable set. The
number of vertices in a maximum stable set, denoted by
$\beta_0(\mathcal{C})$, is
called the {\it stability number} of $\mathcal{C}$. Notice 
that $\alpha_0(\mathcal{C})+\beta_0(\mathcal{C})=n$.

A clutter
$\mathcal{C}$ is called {\it indecomposable\/} if it cannot be decomposed as a
disjoint union of induced subclutters 
$\mathcal{C}_1$, $\mathcal{C}_2$ such that
$\alpha_0(\mathcal{C})=\alpha_0(\mathcal{C}_1)+\alpha_0(\mathcal{C}_2)$
(see Definition \ref{reduced-clutter}). Erd\"os and Gallai
\cite{erdos-gallai} introduced this notion for graphs. 
A clutter obtained from $\mathcal{C}$ by a sequence of
deletions and duplications of vertices is called 
a {\it parallelization} (see Definition~\ref{parallelization-def}). 
If $a=(a_i)$ is a vector in $\mathbb{N}^n$, we denote by $\mathcal{C}^a$
the clutter obtained from
$\mathcal{C}$ by successively deleting any vertex $x_i$ with $a_i=0$ and
duplicating $a_i-1$ times any vertex $x_i$ if $a_i\geq 1$ (see 
Example \ref{april9-09}).  

Our main result is:

\noindent {\bf Theorem~\ref{irr-b-cover-char}}{\it\ Let $0\neq a=(a_i)\in\mathbb{N}^n$,
$b\in\mathbb{N}$. Then $x_1^{a_1}\cdots x_n^{a_n}t^b$ is a minimal
generator of $R_s(I)$, as
a $K$-algebra, if and only if $\mathcal{C}^a$ is an indecomposable clutter
and $b=\alpha_0(\mathcal{C}^a)$.}

There are two cases where a combinatorial description of the symbolic
Rees algebra is known. If the clutter $\mathcal{C}$ has the
max-flow  
min-cut property, then by a result of \cite{clutters}, 
we have $I^i=I^{(i)}$ for all $i\geq 1$, i.e., $R_s(I)=R[It]$ and a 
minimal generator of $R_s(I)$ is either a vertex $x_i$ or an
``edge'' $x_et$ with  $e\in E(\mathcal{C})$. If
$\mathcal{C}$ is a perfect graph, then the minimal generators of
$R_s(I)$ are in one to one correspondence with the cliques
(complete subgraphs) of $\mathcal{C}$ \cite{perfect}.  Both cases 
will follow from our combinatorial description 
of $R_s(I)$ (see Corollaries \ref{coro1-1-irr-b-cover-char} and 
\ref{coro2-irr-b-cover-char} respectively).

As a byproduct, in Section~\ref{irreducible-parallelizations} we
give a method---based on the computation of Hilbert bases of
polyhedral cones---to compute all indecomposable parallelizations of any clutter
$\mathcal{C}$ along with all the corresponding vertex covering numbers. 
In particular our method allows to compute all
indecomposable induced subclutters of any clutter $\mathcal{C}$. 
This means that the symbolic Rees algebra of
$I$ encodes combinatorial information of the clutter
which can be decoded using a computer program, such as {\sc
Normaliz} \cite{normaliz2}, which is able to
compute Hilbert bases of polyhedral cones.

Harary and Plummer \cite{harary-plummer} studied some properties of
indecomposable graphs. They
showed that if a connected graph is separated by the points of a
complete subgraph, then $G$ is decomposable. All indecomposable graphs with
at least three vertices contain at least one odd cycle, and the join
of two 
indecomposable graphs
is indecomposable \cite{harary-plummer}. Indecomposable graphs were first
studied from an algebraic point of view in
\cite{covers}. To the best of
our knowledge there is no structure theorem for indecomposable graphs. 

Indecomposable subgraphs occur naturally in the theory of perfect
graphs. Indeed, a graph $G$ is perfect if and
only if the indecomposable parallelizations
of $G$ are exactly the complete subgraphs or cliques of
$G$ (see Proposition~\ref{perfect-char-parallel}). This was first
shown in \cite{covers} using the main result of \cite{seymour}. For graphs, we can
use our methods to locate all {\it induced odd cycles\/} ({\it odd
holes\/}) and all {\it induced 
complements of odd cycles\/} ({\it odd antiholes}) of length at 
least five. Indeed, odd holes
of any length and odd antiholes of length at 
least five are indecomposable subgraphs (see
Lemma~\ref{antiholes-are-irred}), and thus by
Theorem~\ref{irr-b-cover-char} they correspond to minimal
generators of the symbolic Rees algebra of the edge ideal of the
graph. 
Odd holes and odd antiholes play a 
major role in graph theory. In \cite{seymour} it is shown that a graph $G$ is 
perfect if and only if $G$ is a {\it Berge graph}, i.e., if and only
if $G$ has no odd holes or odd antiholes of length at 
least five. In commutative algebra odd holes occurred for the first
time in \cite{aron-hoyos}, and later in the description
of $I(G)^{\{2\}}$, the join of an edge ideal of a graph $G$ with 
itself \cite{simis-ulrich}. They also occurred in the description of the
associated primes of powers of ideals of vertex covers of graphs \cite{fhv}.

The problem of finding a minimum vertex cover of a graph is a classical
optimization problem in computer science and is a typical example of
an NP-hard problem. From the point of view of 
computational complexity theory, finding all indecomposable
subgraphs of a given graph using our method is a hard problem 
because to apply this method we must know all minimal vertex
covers (see Section~\ref{irreducible-parallelizations}). Thus, 
although our results provide some tools for computing,  
the contributions of this paper could be more interesting from the
theoretical point of view. 

Throughout the paper we introduce most of the 
notions that are relevant for our purposes. For unexplained
terminology, we refer to \cite{diestel,Schr2,Vas,monalg}. 

\section{Symbolic Rees algebras of  edge ideals}\label{blowupcovers}

In this section we will give a combinatorial description of the 
minimal generators of the symbolic Rees 
algebra of the edge ideal of a clutter using the notion 
of a parallelization of a clutter and the notion of an indecomposable
clutter. We continue using the
definitions and terms from the introduction.

Let $\mathcal{C}$ be a clutter with vertex set $X=\{x_1,\ldots,x_n\}$ and let 
$I=I(\mathcal{C})$ be its edge ideal. We denote  by
$\Upsilon(\mathcal{C})$ the clutter whose
edges are the minimal vertex covers of $\mathcal{C}$. The clutter
$\Upsilon(\mathcal{C})$ is called the {\it blocker\/} of
$\mathcal{C}$ or the {\it Alexander dual\/} of $\mathcal{C}$. As
usual, we use  $x^a$ as an abbreviation for $x_1^{a_1} \cdots x_n^{a_n}$,  
where $a=(a_i)\in \mathbb{N}^n$. 

If $C$ is a subset of $X$, its {\it
characteristic vector\/} is the vector $v=\sum_{x_i\in C}e_i$, where
$e_i$ is the $i${\it th} unit vector in $\mathbb{R}^n$. Let 
$C_1,\ldots,C_s$ be the minimal vertex covers of $\mathcal{C}$ and 
let $u_k$ be the characteristic vector of $C_k$ for
$1\leq k\leq s$. In our situation, according to \cite[Proposition~7.3.14]{monalg}, the
$b${\it th} 
symbolic power of $I$ has a simple expression: 
\begin{eqnarray}
I^{(b)}&=&\mathfrak{p}_1^b\cap\cdots\cap\mathfrak{p}_s^b\nonumber\\ 
&=&
(\{x^a\vert\, \langle a,u_k\rangle\geq b\mbox{ for
}k=1,\ldots,s\}),\label{march5-09}
\end{eqnarray}
where $\mathfrak{p}_k$ is the prime ideal of $R$ generated by $C_k$
and $\langle\ ,\, \rangle$ denotes the standard 
inner product. In particular, if $b=1$, we obtain the primary 
decomposition of $I$ because $I^{(1)}=I$. Thus the height of $I$
equals $\alpha_0(\mathcal{C})$, the vertex covering number of
$\mathcal{C}$. This is a hint of the rich interaction between the 
combinatorics of $\mathcal{C}$ and the algebra of $I$.

Next, in Lemma~\ref{april8-09}, we give a simple description of
$R_s(I)$ that was first observed  in the discussion of symbolic Rees algebras given in
\cite[p.~75]{normali}, 
see also \cite{cover-algebras}. 
Let $a=(a_i)\neq 0$ be a 
vector in $\mathbb{N}^n$ and let $b\in\mathbb{N}$. From
Eq.~(\ref{march5-09}) we get that $x^at^b$ is in $R_s(I)$ if and only
if 
$$
\langle a,u_k\rangle\geq b\ \mbox{for}\ k=1,\ldots,s.
$$  
If $a, b$ satisfy this system of linear inequalities, we say that 
$a$ is a $b$-{\it vertex cover\/} of
$\Upsilon(\mathcal{C})$. Often we will call a $b$-vertex cover simply a $b$-{\it
cover\/}.  Thus the symbolic Rees algebra of $I$ is equal to the 
$K$-subalgebra of $R[t]$ generated by all monomials $x^at^b$ such 
that $a$ is a $b$-cover of $\Upsilon(\mathcal{C})$, as was first
shown in \cite[Theorem~3.5]{normali}. The notion of a $b$-cover 
occurs in combinatorial optimization (see for instance
\cite[Chapter~77, p.~1378]{Schr2} and the references 
there) and algebraic 
combinatorics \cite{normali,cover-algebras}. We say that a
$b$-cover $a$ of $\Upsilon(\mathcal{C})$ is
{\it decomposable\/} if there exists an $i$-cover $c$ and a $j$-cover
$d$ of $\Upsilon(\mathcal{C})$ 
such that $a=c+d$ and $b=i+j$. If $a$ is not decomposable, we call $a$
{\it indecomposable\/}. The indecomposable $0$ and $1$ 
covers of $\Upsilon(\mathcal{C})$ are the unit vectors $e_1,\ldots,e_n$ and 
the characteristic vectors $v_1,\ldots,v_q$ of the edges of $\mathcal{C}$
respectively. 

\begin{lemma}\label{april8-09} A monomial $x^at^b\neq 1$ is a 
minimal generator of $R_s(I)$, as a $K$-algebra, if and only if $a$ is an indecomposable
$b$-cover of $\Upsilon(\mathcal{C})$. In particular, the following
equality holds{\rm :}
\begin{equation}\label{sept5-07}
R_s(I)=K[\{x^{a}t^b\vert\, a\mbox{ is an indecomposable }b\mbox{-cover
of }  \Upsilon(\mathcal{C})\}].
\end{equation}
\end{lemma}

\begin{proof} It follows from the discussion above, by
decomposing any $b$-cover into indecomposable ones.
\end{proof}

Let $S$ be a set of vertices of a clutter $\mathcal{C}$.  The {\it induced
subclutter\/} on $S$, denoted by $\mathcal{C}[S]$, is the maximal
subclutter of $\mathcal{C}$ with vertex set $S$. Thus the vertex set
of $\mathcal{C}[S]$ is $S$ and the edges of $\mathcal{C}[S]$ are
exactly the edges of $\mathcal{C}$ contained in $S$. Notice that
$\mathcal{C}[S]$ may have isolated vertices, i.e., vertices that do
not belong to any edge of $\mathcal{C}[S]$. If $\mathcal{C}$ is a
discrete clutter, i.e., 
all the vertices of $\mathcal{C}$ are
isolated, we set $I(\mathcal{C})=0$ and $\alpha_0(\mathcal{C})=0$.

Let $\mathcal{C}$ be a clutter and let $X_1,X_2$ be a partition of
$V(\mathcal{C})$ into nonempty sets. Clearly, one has the inequality 
\begin{equation}\label{may1-09}
\alpha_0(\mathcal{C})\geq \alpha_0(\mathcal{C}[X_1])+\alpha_0(\mathcal{C}[X_2]). 
\end{equation}
If $\mathcal{C}$ is a graph and equality occurs, Erd\"os and 
Gallai \cite{erdos-gallai} call $\mathcal{C}$ a {\it decomposable
graph\/}. This motivates the following similar notion for clutters. 
\begin{definition}\label{reduced-clutter}\rm A clutter $\mathcal{C}$ is called {\it
decomposable} if there are nonempty vertex sets $X_1, X_2$ such that 
$X$ is the disjoint union of $X_1$ and $X_2$, 
and
$\alpha_0(\mathcal{C})=\alpha_0(\mathcal{C}[X_1])+
\alpha_0(\mathcal{C}[X_2])$. If $\mathcal{C}$ is not decomposable, it is
called {\it indecomposable\/}. 
\end{definition}

Examples of indecomposable graphs include complete graphs, odd cycles 
and complements of odd cycles of length at least five 
(see Lemma~\ref{antiholes-are-irred}).

\begin{definition}{(Schrijver \cite{Schr2})}\label{parallelization-def}
The {\it 
duplication\/} of a vertex $x_i$  
of a clutter $\mathcal{C}$ means extending its vertex set $X$ by a
new vertex $x_i'$ 
and replacing
$E(\mathcal{C})$ by
$$    
E(\mathcal{C})\cup\{(e\setminus\{x_i\})\cup\{x_i'\}\vert\, x_i\in e\in
E(\mathcal{C})\}.
$$
The {\it deletion\/} of $x_i$, denoted by
$\mathcal{C}\setminus\{x_i\}$,  is the clutter formed from
$\mathcal{C}$ by deleting the vertex $x_i$ and all edges containing
$x_i$. A clutter obtained from $\mathcal{C}$ by a sequence of
deletions and 
duplications of vertices is called a {\it parallelization\/}. 
\end{definition} 

It is not difficult to verify that these two operations commute. 
If $a=(a_i)$ is a vector in $\mathbb{N}^n$, we denote by $\mathcal{C}^a$
the clutter obtained from
$\mathcal{C}$ by successively deleting any vertex $x_i$ with $a_i=0$ and
duplicating $a_i-1$ times any vertex $x_i$ if $a_i\geq 1$ (for graphs
cf. \cite[p.~53]{golumbic}).

\begin{example}\label{april9-09} Let $G$ be the graph whose only edge is 
$\{x_1,x_2\}$ and let $a=(3,3)$. We set $x_i^1=x_i$ for $i=1,2$. The parallelization 
$G^a$ is a complete bipartite graph 
with bipartition $V_1=\{x_1^1,x_1^2,x_1^3\}$ and
$V_2=\{x_2^1,x_2^2,x_2^3\}$. Note that $x_i^k$ is a vertex, i.e., $k$
is an index not an exponent.
\vspace{0.8cm}
$$
\begin{array}{cccc}
\setlength{\unitlength}{.04cm}
\thicklines
\begin{picture}(80,35)
\put(10,10){\circle*{3.1}}
\put(10,40){\circle*{3.1}}
\put(-5,42){$x_1$}
\put(-5,3){$x_2$}
\put(10,10){\line(0,1){30}}
\put(20,15){$G$}
\put(-10,-15){\mbox{Fig. 1. Graph}}
\end{picture}
&
\setlength{\unitlength}{.04cm}
\thicklines
\begin{picture}(80,35)
\put(30,40){\circle*{3.1}}
\put(60,40){\circle*{3.1}}
\put(0,40){\circle*{3.1}}
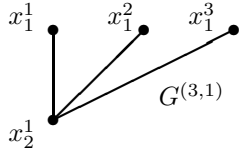
\put(-20,-15){\mbox{Fig. 2. Duplications of $x_1$}}
\put(0,10){\circle*{3.1}}
\put(-15,42){$x_1^1$}
\put(18,42){$x_1^2$}
\put(45,42){$x_1^3$}
\put(-15,3){$x_2^1$}
\put(0,10){\line(0,1){30}}
\put(0,10){\line(1,1){30}}
\put(0,10){\line(2,1){60}}
\put(35,15){$G^{(3,1)}$}
\end{picture}
&
\ \ \ \ \ \ \ \ \ \ 
&
\setlength{\unitlength}{.04cm}
\thicklines
\begin{picture}(80,35)
\put(0,10){\circle*{3.1}}
\put(-15,42){$x_1^1$}
\put(18,42){$x_1^2$}
\put(45,42){$x_1^3$}
\put(-15,3){$x_2^1$}
\put(18,3){$x_2^2$}
\put(45,3){$x_2^3$}
\put(0,10){\line(0,1){30}}
\put(0,10){\line(1,1){30}}
\put(0,10){\line(2,1){60}}
\put(30,10){\circle*{3.1}}
\put(30,10){\line(0,1){30}}
\put(30,10){\line(1,1){30}}
\put(30,10){\line(-1,1){30}}
\put(60,10){\circle*{3.1}}
\put(60,10){\line(0,1){30}}
\put(60,10){\line(-2,1){60}}
\put(60,10){\line(-1,1){30}}
\put(0,40){\circle*{3.1}}
\put(30,40){\circle*{3.1}}
\put(60,40){\circle*{3.1}}
\put(70,15){$G^{(3,3)}$}
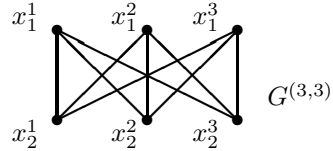
\put(-15,-15){\mbox{Fig. 3. Duplications of $x_2$}}
\end{picture}
\end{array}
$$
\end{example}
\vspace{0.6cm}

\begin{proposition}{(\cite[Lemma~2.15]{cm-mfmc},
\cite[p.~1385, Eq.~(78.6)]{Schr2})}\label{alpha-exp} 
Let $\mathcal{C}$ be a clutter with $n$ vertices and let  
$\Upsilon(\mathcal{C})$ be the blocker of 
$\mathcal{C}$. If $a=(a_i)\in\mathbb{N}^n$, then
$$
\left.\min\left\{\sum_{x_i\in
C}a_i\right\vert\,
C\in\Upsilon(\mathcal{C})\right\}
= \alpha_0(\mathcal{C}^a).
$$
\end{proposition}

We come to the main result of this section.

\begin{theorem}\label{irr-b-cover-char} Let $\mathcal{C}$ be a clutter with vertex set
$X=\{x_1,\ldots,x_n\}$ and let $0\neq a=(a_i)\in\mathbb{N}^n$,
$b\in\mathbb{N}$. Then $x^at^b$ is a minimal generator of
$R_s(I(\mathcal{C}))$, as a $K$-algebra, if and only if\/
$\mathcal{C}^a$ is an indecomposable clutter and 
$b=\alpha_0(\mathcal{C}^a)$.  
\end{theorem}

\begin{proof} 
We may assume that 
$a=(a_1,\ldots,a_m,0,\ldots,0)$, where $a_i\geq 1$ for 
$i=1,\ldots,m$. For each $1\leq i\leq m$ the vertex $x_i$ is 
duplicated $a_i-1$ times, and the vertex $x_i$ is deleted for
each $i>m$. 
We denote the duplications of $x_i$ by 
$x_{i}^2,\ldots,x_i^{a_i}$ and set $x_i^1=x_i$ for $1\leq i\leq m$.
Thus the vertex set of  
$\mathcal{C}^a$ can be written as 
$$
X^a=\{x_1^1,\ldots,x_1^{a_1},\ldots,x_i^1,\ldots,x_i^{a_i},
\ldots,x_m^1,\ldots,x_m^{a_m}\}=X^{a_1}
\cup X^{a_2}\cup\cdots\cup X^{a_m},
$$
where $X^{a_i}=\{x_i^1,\ldots,x_i^{a_i}\}$ for $1\leq i\leq m$ and 
$X^{a_i}\cap X^{a_j}=\emptyset$ for $i\neq j$. 

$\Rightarrow$) Assume that $x^at^b$ is a minimal generator of
$R_s(I(\mathcal{C}))$. Then, by Lemma~\ref{april8-09}, 
$a$ is an indecomposable $b$-cover of
$\Upsilon(\mathcal{C})$. First we prove that
$b=\alpha_0(\mathcal{C}^a)$. There is $k$ such that $a_k\neq 0$. We
may assume that $a-e_k\neq 0$. By Proposition~\ref{alpha-exp} we need
only show the equality 
$$
b=\left.\min\left\{\sum_{x_i\in
C}a_i\right\vert\,
C\in\Upsilon(\mathcal{C})\right\}.
$$
As $a$ is a $b$-cover of $\Upsilon(\mathcal{C})$, the minimum is
greater than or equal to $b$. If the minimum is greater than $b$, then we
can write $a=(a-e_k)+e_k$, where $a-e_k$ is a $b$-cover and $e_k$ is a
$0$-cover, a contradiction to the indecomposability of $a$. 

Next we show
that $\mathcal{C}^a$ is indecomposable. 
We proceed by contradiction.
Assume that $\mathcal{C}^a$ is decomposable. Then there is a
partition $X_1, X_2$ of $X^a$ such that 
$\alpha_0(\mathcal{C}^a)=\alpha_0(\mathcal{C}^a[X_1])+
\alpha_0(\mathcal{C}^a[X_2])$. For $1\leq i\leq n$, we set
$$
\ell_i=|X^{a_i}\cap X_1|\ \mbox{ and }\ p_i=|X^{a_i}\cap X_2|
$$
if $1\leq i\leq m$ and $\ell_i=p_i=0$ if $i>m$. Consider the 
vectors $\ell=(\ell_i)$ and $p=(p_i)$. Notice that $a$ has a 
decomposition $a=\ell+p$ because one has a partition
$X^{a_i}=(X^{a_i}\cap X_1)\cup(X^{a_i}\cap X_2)$ for $1\leq i\leq m$.
To derive a contradiction we now claim that
$\ell$ (resp. $p$) is an $\alpha_0(\mathcal{C}^a[X_1])$-cover 
(resp. $\alpha_0(\mathcal{C}^a[X_2])$-cover) of $\Upsilon(\mathcal{C})$.
Take an arbitrary $C$ in $\Upsilon(\mathcal{C})$. The set
$$
C_a=\bigcup_{x_i\in C}\{x_i^1,\ldots,x_i^{a_i}\}=\bigcup_{x_i\in C}X^{a_i}
$$
is a vertex cover of $\mathcal{C}^a$. Indeed, if $f_k$ is any edge of
$\mathcal{C}^a$, then $f_k$ has the form
\begin{equation}\label{may6-09}
f_k=\{x_{k_1}^{j_{k_1}},x_{k_2}^{j_{k_2}},\ldots,x_{k_r}^{j_{k_r}}\}\ \
\ \ \ 
(1\leq k_1<\cdots<k_r\leq m;\ 1\leq j_{k_i}\leq a_{k_i})
\end{equation}
for some edge $\{x_{k_1},x_{k_2},\ldots,x_{k_r}\}$ of $\mathcal{C}$.
Since $\{x_{k_1},x_{k_2},\ldots,x_{k_r}\}\cap C\neq\emptyset$, 
we get $f_k\cap C_a\neq\emptyset$. Thus $C_a$ is a vertex cover of
$\mathcal{C}^a$. Therefore $C_a\cap X_1$ and $C_a\cap X_2$ are vertex
covers of $\mathcal{C}^a[X_1]$ and $\mathcal{C}^a[X_2]$ respectively
because $E(\mathcal{C}^a[X_i])$ is contained in $E(\mathcal{C}^a)$ for
$i=1,2$. 
Hence using the partitions
$$
C_a\cap X_1=\bigcup_{x_i\in C}(X^{a_i}\cap X_1)\ \ \mbox{ and }\ \ 
C_a\cap X_2=\bigcup_{x_i\in C}(X^{a_i}\cap X_2)
$$ 
we obtain
$$
\alpha_0(\mathcal{C}^a[X_1])
\leq|C_a\cap X_1|=
\sum_{x_i\in
C}\ell_i 
\ \ 
\mbox{ and }\ \ \alpha_0(\mathcal{C}^a[X_2])\leq|C_a\cap X_2|=\sum_{x_i\in
C}p_i.
$$
This completes the proof of the claim. Consequently $a$ is
a decomposable $b$-cover of $\Upsilon(\mathcal{C})$, where
$b=\alpha_0(\mathcal{C}^a)$, a contradiction to
the indecomposability of $a$. 

$\Leftarrow$) Assume that $\mathcal{C}^a$ is an indecomposable clutter
and $b=\alpha_0(\mathcal{C}^a)$. To show that $x^at^b$ is a mi\-nimal
generator of $R_s(I(\mathcal{C}))$ we need only show that $a$ is an 
indecomposable $b$-cover of
$\Upsilon(\mathcal{C})$.  
To begin with, notice that $a$ is a
$b$-cover of $\Upsilon(\mathcal{C})$ by Proposition~\ref{alpha-exp}. 
We proceed by contradiction assuming that there is a
decomposition $a=\ell+p$, where $\ell=(\ell_i)$ is a $c$-cover of
$\Upsilon(\mathcal{C})$, $p=(p_i)$ is a $d$-cover of
$\Upsilon(\mathcal{C})$, and $b=c+d$. Each 
$X^{a_i}$ can be decomposed as $X^{a_i}=X^{\ell_i}\cup X^{p_i}$, 
where $X^{\ell_i}\cap X^{p_i}=\emptyset$, $\ell_i=|X^{\ell_i}|$, 
and  $p_i=|X^{p_i}|$. We set 
$$
X^{\ell}=X^{\ell_1}\cup\cdots\cup
X^{\ell_m}\ \  \mbox{ and }\ \ X^{p}=X^{p_1}\cup\cdots\cup X^{p_m}. 
$$
Then one has a decomposition $X^a=X^\ell\cup X^p$ of the vertex set of
$\mathcal{C}^a$. We now show that $\alpha_0(\mathcal{C}^a[X^\ell])\geq c$ 
and $\alpha_0(\mathcal{C}^a[X^p])\geq d$. By symmetry, it suffices to
prove the first inequality. Take an arbitrary minimal vertex cover
$C_\ell$ of $\mathcal{C}^a[X^\ell]$. Then $C_\ell\cup´X^p$ is a vertex
cover of $\mathcal{C}^a$ because if $f$ is an edge of $\mathcal{C}^a$
contained in $X^\ell$, then $f$ is covered by $C_\ell$, otherwise $f$
is covered by $X^p$. Hence there is a minimal vertex cover $C_a$ of 
$\mathcal{C}^a$ such that $C_a\subset C_\ell\cup X^p$. Since
$\mathcal{C}[\{x_1,\ldots,x_m\}]$ is a subclutter of $\mathcal{C}^a$,
there is a minimal vertex cover $C_1$ of
$\mathcal{C}[\{x_1,\ldots,x_m\}]$ contained in $C_a$. 
Then the set $C_1\cup\{x_i\vert\, i>m\}$ is a vertex cover of
$\mathcal{C}$. Therefore there is a minimal vertex cover $C$ of
$\mathcal{C}$ such that $C\cap\{x_1,\ldots,x_m\}\subset C_a$.
Altogether one has: 
\begin{eqnarray}
&&C\cap\{x_1,\ldots,x_m\}\subset C_a\subset C_\ell\cup X^p\Longrightarrow\label{mar4-09}\\
&&\ \ \ \ \ \ C\cap\{x_1,\ldots,x_m\}
\subset C_a\cap\{x_1,\ldots,x_m\}\subset (C_\ell\cup X^p)\cap
\{x_1,\ldots,x_m\}.
\end{eqnarray}
We may assume that $C_a\cap\{x_1,\ldots,x_m\}=\{x_1,\ldots,x_s\}$.
Next we claim that $X^{a_i}\subset C_a$ for $1\leq i\leq s$. Take an
integer $i$ between $1$ and $s$. Since $C_a$ is a minimal vertex
cover of 
$\mathcal{C}^a$, there exists
an edge $e$ of $\mathcal{C}^a$ such that $e\cap C_a=\{x_i^1\}$. Then
$(e\setminus\{x_i^1\})\cup\{x_i^j\}$ is an edge of $\mathcal{C}^a$ for
$j=1,\ldots, a_i$, this follows using that the edges of $\mathcal{C}^a$
are of the form described in Eq.~(\ref{may6-09}). Consequently
$x_i^j\in C_a$ for $j=1,\ldots,a_i$.  
This completes the proof of the claim. Thus one has 
$X^{\ell_i}\subset X^{a_i}\subset C_a$ for $1\leq i\leq s$. Hence, 
by Eq.~(\ref{mar4-09}), and noticing that $X^{\ell_i}\cap X^p=\emptyset$, 
we get
$X^{\ell_i}\subset C_\ell$ for $1\leq i\leq s$. So, using that 
$\ell_i=0$ for $i>m$, we get
$$
\alpha_0(\mathcal{C}^a[X^\ell])\geq|C_\ell|\geq\sum_{i=1}^s\ell_i\geq \sum_{x_i\in
C\cap\{x_1,\ldots,x_m\}}\hspace{-5mm}\ell_i=\sum_{x_i\in C}\ell_i\geq
c.
$$
Therefore $\alpha_0(\mathcal{C}^a[X^\ell])\geq c$. Similarly
$\alpha_0(\mathcal{C}^a[X^p])\geq d$. Thus 
$$
\alpha_0(\mathcal{C}^a[X^\ell])+\alpha_0(\mathcal{C}^a[X^p])\geq
c+d=b,
$$ 
and consequently by Eq.~(\ref{may1-09}) we have the equality
$$
\alpha_0(\mathcal{C}^a[X^\ell])+\alpha_0(\mathcal{C}^a[X^p])=
\alpha_0(\mathcal{C}^a).
$$
Thus we have shown that $\mathcal{C}^a$ is a decomposable
clutter, a contradiction.
\end{proof}

Let $\mathcal{C}$ be a clutter. A set of edges of $\mathcal C$ is
called {\it independent\/} if no two of them have a
common vertex. We denote the maximum number of independent edges of ${\mathcal C}$ by 
$\beta_1({\mathcal C})$, this number is called the {\it matching
number\/} of $\mathcal{C}$. In general the vertex covering number and the
matching number satisfy $\beta_1(\mathcal{C})\leq
\alpha_0(\mathcal{C})$.

\begin{definition}\label{konig-def}\rm If $\beta_1({\mathcal
C})=\alpha_0({\mathcal C})$, we say  that $\mathcal C$ has the {\it K\"onig
property\/}.
\end{definition}

\begin{lemma}\label{konig+irred=edge} If $\mathcal{C}$ is an
indecomposable clutter with the K\"onig property, then 
either $\mathcal{C}$ has no edges and has exactly one isolated vertex or
$\mathcal{C}$ has only one edge and no isolated vertices. 
\end{lemma}

\begin{proof} Let $f_1,\ldots,f_g$ be a set of independent edges and
let $X'=\cup_{i=1}^gf_i$, where 
$g=\alpha_0(\mathcal{C})$. Note that $g=0$ if $\mathcal{C}$ has no
edges. Then $V(\mathcal{C})$ has a partition
\[
V(\mathcal{C})=\left(\displaystyle\cup_{i=1}^gf_i\right)\cup\left(\cup_{x_i\in
V(\mathcal{C})\setminus X'}\{x_i\}\right).
\]
As $\mathcal{C}$ is indecomposable, we get that either $g=0$ and 
$V(\mathcal{C})=\{x_i\}$ for some vertex $x_i$ or $g=1$ and 
$V(\mathcal{C})=f_i$ for some $i$. Thus in the second case, as $\mathcal{C}$ is a
clutter, we get that $\mathcal{C}$ has exactly one edge and no
isolated vertices. 
\end{proof}

\begin{corollary}\label{coro1-irr-b-cover-char} 
Let $\mathcal{C}$ be a clutter and let
$I=I(\mathcal{C})$ be its edge ideal. Then all indecomposable 
para\-llelizations of $\mathcal{C}$ satisfy the K\"onig property if
and only if $I^i=I^{(i)}$ for $i\geq 1$.
\end{corollary}

\begin{proof} $\Rightarrow$) It suffices to prove that
$R[It]=R_s(I)$. 
Clearly 
$R[It]\subset R_s(I)$. To prove the reverse inclusion take a minimal
generator $x^at^b$ of $R_s(I)$. If $b=0$, then $a=e_i$ for some $i$
and $x^at^b=x_i$. Thus $x^at^b\in R[It]$. Assume $b\geq 1$. 
By Theorem~\ref{irr-b-cover-char} we
have that $\mathcal{C}^a$ is an indecomposable clutter such that 
$b=\alpha_0(\mathcal{C}^a)$. As $\mathcal{C}^a$ is indecomposable and 
satisfies the K\"onig property, using Lemma~\ref{konig+irred=edge},
it is not hard to see that $b=1$ and 
that $E(\mathcal{C}^a)=\{e\}$ consists of a single edge $e$ of
$\mathcal{C}$, i.e., $x^at^b=x_et$, where $x_e=\prod_{x_i\in e}x_i$. 
Thus $x^at^b\in R[It]$. 

$\Leftarrow$) Since $R[It]=R_s(I)$, by Theorem~\ref{irr-b-cover-char}
we obtain that the only indecomposable parallelizations are either induced 
subclutters of $\mathcal{C}$ with exactly one edge and no isolated
vertices or subclutters consisting of exactly one isolated vertex. 
Thus in both cases they satisfy the K\"onig property.
\end{proof}

A clutter $\mathcal{C}$ is called {\it Mengerian\/} if all its
parallelizations have the K\"onig property. A clutter $\mathcal{C}$ satisfies 
the  {\em max-flow min-cut property} if the linear program:
\begin{equation*}
\textup{max} \{ \langle\mathbf{1},y\rangle \, \vert \, y \geq 0, \, A y \leq a \}
\end{equation*} 
has an integral optimal solution for all $a \in \mathbb{N}^n$, where $A$ is
the incidence matrix of the clutter $\mathcal{C}$ and $\mathbf{1}$ is the vector
of all ones. The columns of $A$ are the
characteristic vectors of the edges of $\mathcal{C}$. It is well known
that a clutter is Mengerian if and only if it satisfies the max-flow
min-cut property \cite[Chapter 79]{Schr2}.

Thus the last corollary can be restated as:

\begin{corollary}{\rm \cite[Corollary~3.14]{clutters}}\label{coro1-1-irr-b-cover-char} 
Let $\mathcal{C}$ be a clutter and let
$I$ be its edge ideal. Then $\mathcal{C}$ has the max-flow min-cut
property if and only if $I^i=I^{(i)}$ for $i\geq 1$.
\end{corollary}

The following was the first deep result in the study of 
symbolic powers of edge ideals from the viewpoint of graph theory. 

\begin{corollary}{\rm \cite[Theorem
5.9]{ITG}}\label{coro3-irr-b-cover-char}  
Let $G$ be a graph and let
$I$ be its edge ideal. Then $G$ is bipartite if and only if 
$I^i=I^{(i)}$ for $i\geq 1$.
\end{corollary}

\begin{proof} $\Rightarrow$) If $G$ is a bipartite graph, then any
parallelization of $G$ is again a bipartite graph. This means that any
parallelization of $G$ satisfies the K\"onig property because 
bipartite graphs satisfy this property \cite[Theorem 2.1.1]{diestel}. 
Thus $I^i=I^{(i)}$ for all $i$ 
by Corollary \ref{coro1-irr-b-cover-char}.

$\Leftarrow$) Assume that $I^i=I^{(i)}$ for $i\geq 1$. By 
Corollary \ref{coro1-irr-b-cover-char} all indecomposable induced 
subgraphs of $G$ have the K\"onig property. If $G$ is not
bipartite, then $G$ has an induced odd cycle, a contradiction 
because induced odd cycles are indecomposable \cite{harary-plummer} 
and do not satisfy the
K\"onig property. 
\end{proof}

\begin{corollary}\label{march20-09} Let $\mathcal{C}$ be a clutter with
vertex set $X=\{x_1,\ldots,x_n\}$ and let $S\subset X$. Then the
induced clutter $H=\mathcal{C}[S]$ is indecomposable if and only if 
the monomial $(\prod_{x_i\in S}x_i)t^{\alpha_0(H)}$ is a minimal
generator of $R_s(I(\mathcal{C}))$.
\end{corollary} 

\begin{proof} Let $a=\sum_{x_i\in S}e_i$. Since
$\mathcal{C}^a=\mathcal{C}[S]$, the result follows from
Theorem~\ref{irr-b-cover-char}.
\end{proof}

\begin{corollary} Let $\mathcal{C}$ be a clutter with $n$ vertices and
let $A$ be its incidence matrix. If the polyhedron 
$Q(A)=\{x\vert\, x\geq
0;xA\geq\mathbf{1}\}$ has only integral vertices,  
then $\alpha_0(\mathcal{C}^a)\leq n-1$ for all indecomposable 
parallelizations $\mathcal{C}^a$ of $\mathcal{C}$.
\end{corollary}

\begin{proof} Let $v_1,\ldots,v_q$ be the characteristic vectors of 
the edges of $\mathcal{C}$ and let $\overline{I^i}$ be the 
integral closure of $I^i$, where $I$ is the edge ideal of $\mathcal{C}$. 
As $Q(A)$ is integral, by \cite[Corollary~3.13]{clutters} we have that 
$\overline{I^i}=I^{(i)}$ for $i\geq 1$, where 
$$\overline{I^i}=(\{x^a\vert\, 
a\in{iB}\cap\mathbb{Z}^n\})
$$
and $B=\mathbb{Q}_+^n+{\rm conv}(v_1,\ldots,v_q)$, see \cite{monalg}.
Thus we have the
equality $\overline{R[It]}=R_s(I)$, where $\overline{R[It]}$ is the
integral closure of $R[It]$ in its field of fractions. 
Take any indecomposable parallelization $\mathcal{C}^a$ of $\mathcal{C}$ and
consider the monomial $m=x^at^b$, where $b=\alpha_0(\mathcal{C}^a)$. By
Theorem~\ref{irr-b-cover-char} $m$ is a minimal generator of 
$R_s(I)$. Now, according to \cite[Corollary~3.11]{normali}, a minimal generator 
of $\overline{R[It]}$ has degree in $t$ at most $n-1$, i.e., 
$b\leq n-1$.        
\end{proof}

We end this section showing some very basic properties of indecomposable
clutters. If $e$ is a edge of a clutter $\mathcal{C}$, we
denote by $\mathcal{C}\setminus\{e\}$ the spanning subclutter of
$\mathcal{C}$ obtained by
deleting $e$ and keeping all the vertices of $\mathcal{C}$.

\begin{definition}\rm A clutter $\mathcal{C}$ is called {\it vertex
critical\/} if 
$\alpha_0{(\mathcal{C}\setminus\{x_i\})}<\alpha_0{(\mathcal{C})}$ 
for all $x_i\in V(\mathcal{C})$. A clutter $\mathcal{C}$ is 
called {\it edge critical\/} if 
$\alpha_0(\mathcal{C}\setminus\{e\})<\alpha_0{(\mathcal{C})}$ for all
$e\in E(\mathcal{C})$.         
\end{definition}

The next lemma is not hard to prove. 

\begin{lemma}\label{nov29-07}\rm Let $x_i$ be a vertex of a clutter
$\mathcal{C}$ and let $e$ be an edge of $\mathcal{C}$. 

(a) If 
$\alpha_0{(\mathcal{C}\setminus\{x_i\})}<\alpha_0{(\mathcal{C})}$, then
$\alpha_0(\mathcal{C}\setminus\{x_i\})=\alpha_0(\mathcal{C})-1$. 

(b)
If $\alpha_0{(\mathcal{C}\setminus\{e\})}<\alpha_0{(\mathcal{C})}$, then
$\alpha_0(\mathcal{C}\setminus\{e\})=\alpha_0(\mathcal{C})-1$.  
\end{lemma}

\begin{definition}\rm A clutter $\mathcal{C}$ is called {\it
connected\/} if there is no $U\subset V(\mathcal{C})$ such that 
$\emptyset\subsetneq U\subsetneq V(\mathcal{C})$ and such that 
$e\subset U$ or $e\subset V(\mathcal{C})\setminus U$ for each 
edge $e$ of $\mathcal{C}$. 
\end{definition}

\begin{proposition}\label{march20-09-1} 
If a clutter $\mathcal{C}$ is indecomposable, then it is connected 
and vertex critical.  
\end{proposition}

\begin{proof} Assume that $\mathcal{C}$ is disconnected. Then 
there is a partition $X_1,X_2$ of $V(\mathcal{C})$ such that 
\begin{equation}\label{march21-09}
E(\mathcal{C})\subset E(\mathcal{C}[X_1])\cup E(\mathcal{C}[X_2]).
\end{equation}
For $i=1,2$, let $C_i$ be a minimal vertex cover of $\mathcal{C}[X_i]$
with $\alpha_0(\mathcal{C}[X_i])$ vertices. Then, by
Eq.~(\ref{march21-09}), $C_1\cup C_2$ is a
minimal vertex cover of $\mathcal{C}$. Hence 
$\alpha_0(\mathcal{C}[X_1])+\alpha_0(\mathcal{C}[X_2])$ is greater
than or equal to $\alpha_0(\mathcal{C})$. So $\alpha_0(\mathcal{C})$ is equal to
$\alpha_0(\mathcal{C}[X_1])+\alpha_0(\mathcal{C}[X_2])$, a
contradiction to the indecomposability of $\mathcal{C}$. Thus
$\mathcal{C}$ is connected. 

We now show that 
$\alpha_0(\mathcal{C}\setminus\{x_i\})<\alpha_0(\mathcal{C})$ for all $i$. If
$\alpha_0(\mathcal{C}\setminus\{x_i\})=\alpha_0(\mathcal{C})$, then
$V(\mathcal{C})=X_1\cup X_2$, where
$X_1=V(\mathcal{C})\setminus\{x_i\}$ and 
$X_2=\{x_i\}$. Note that $\mathcal{C}[X_1]=\mathcal{C}\setminus\{x_i\}$. 
As $\alpha_0(\mathcal{C}[X_1])=\alpha_0(\mathcal{C})$
and $\alpha_0(\mathcal{C}[X_2])=0$, we contradict the indecomposability
of $\mathcal{C}$. Thus
$\alpha_0(\mathcal{C}\setminus\{x_i\})<\alpha_0(\mathcal{C})$ and
$\mathcal{C}$ is vertex critical. \end{proof}

\begin{proposition}\label{march20-09-2} If $\mathcal{C}$ is a
connected edge critical clutter, then $\mathcal{C}$ is indecomposable.
\end{proposition}

\begin{proof} Assume that $\mathcal{C}$ is decomposable. Then there is a
partition $X_1,X_2$ of $V(\mathcal{C})$ into nonempty vertex sets such
that $\alpha_0(\mathcal{C})=
\alpha_0(\mathcal{C}[X_1])+\alpha_0(\mathcal{C}[X_2])$. Since
$\mathcal{C}$ is connected, there is an edge $e\in E(\mathcal{C})$
intersecting both $X_1$ and $X_2$. Pick a
minimal vertex cover $C$ of $\mathcal{C}\setminus\{e\}$ with less than
$\alpha_0(\mathcal{C})$ 
vertices. As $E(\mathcal{C}[X_i])$ is a 
subset of $E(\mathcal{C}\setminus\{e\})=E(\mathcal{C})\setminus\{e\}$
for $i=1,2$, 
we get that $C$ covers all edges of $\mathcal{C}[X_i]$ for $i=1,2$. Hence $C$ must have at 
least $\alpha_0(\mathcal{C})$ vertices, a contradiction. \end{proof}

From Propositions~\ref{march20-09-1} and \ref{march20-09-2} we obtain:

\begin{corollary} The following hold for any connected
clutter {\rm :}
$$
\begin{array}{cccccc}
\mbox{edge critical}&\Longrightarrow& \mbox{indecomposable}&
&\Longrightarrow
&\mbox{vertex critical.}
\end{array}
$$
\end{corollary}

The next result can be used to build indecomposable clutters. 

\begin{proposition}\label{building-lemma} Let $\mathcal{D}$ be a
clutter obtained from a clutter $\mathcal{C}$ by adding a new
vertex $v$ and some new edges containing $v$ and some vertices of
$V(\mathcal{C})$. 
If $a=(1,\ldots,1)\in\mathbb{N}^n$ is an indecomposable
$\alpha_0(\mathcal{C})$-cover  of $\Upsilon(\mathcal{C})$ such that
$\alpha_0(\mathcal{D})=\alpha_0(\mathcal{C})+1$, then $a'=(a,1)$
is an indecomposable $\alpha_0(\mathcal{D})$-cover of
$\Upsilon(\mathcal{D})$. 
\end{proposition}

\begin{proof}Clearly $a'$ is an $\alpha_0(\mathcal{D})$-cover of $\Upsilon(\mathcal{D})$. 
Assume that $a'=a_1'+a_2'$, where $a_i'\neq 0$ is a $b_i'$-cover of
$\Upsilon(\mathcal{D})$ and $b_1'+b_2'=\alpha_0(\mathcal{D})$. We may assume that
$a_1'=(1,\ldots,1,0,\ldots,0)$ and $a_2'=(0,\ldots,0,1,\ldots,1)$. Let
$a_i$ be the vector in $\mathbb{N}^n$ obtained from $a_i'$ by 
removing its last entry. Set $v=x_{n+1}$. Take a minimal vertex cover
$C_k$ of $\mathcal{C}$ 
and consider $C_k'=C_k\cup\{x_{n+1}\}$. Let $u_k'$ (resp. $u_k$) be 
the characteristic vector of $C_k'$ (resp. $C_k$). Then
$$
\langle a_1,u_k\rangle=\langle a_1',u_k'\rangle\geq b_1'\ \mbox{and}\ 
\langle a_2,u_k\rangle+1=\langle a_2',u_k'\rangle\geq b_2',
$$
and consequently $a_1$ is a $b_1'$-cover of $\Upsilon(\mathcal{C})$. If
$b_2'=0$, then $a_1$ is an $\alpha_0(\mathcal{D})$-cover of $\Upsilon(\mathcal{C})$, 
a contradiction; because if $u$ is the characteristic vector of a minimal
vertex cover of $\mathcal{C}$ with $\alpha_0(\mathcal{C})$ elements, then we would obtain
$\alpha_0(\mathcal{C})\geq\langle u,a_1\rangle\geq\alpha_0(\mathcal{D})$, which is
impossible. Thus $b_2'\geq 1$, and $a_2$ is a $(b_2'-1)$-cover of
$\Upsilon(\mathcal{C})$ if $a_2\neq 0$.
Hence $a_2=0$, because $a=a_1+a_2$ and $a$ is 
indecomposable. This means that $a_2'=e_{n+1}$ is a $b_2'$-cover of
$\Upsilon(\mathcal{D})$, a contradiction. Therefore $a'$ is an indecomposable
$\alpha_0(\mathcal{D})$-cover of $\Upsilon(\mathcal{D})$, as required. \end{proof}

\section{Indecomposable parallelizations and Hilbert
bases}\label{irreducible-parallelizations}

Let $\mathcal{C}$ be a clutter with vertex set $X=\{x_1,\ldots,x_n\}$
and let $C_1,\ldots,C_s$ be the minimal vertex covers of
$\mathcal{C}$. For $1\leq k\leq n$, we denote the characteristic
vector of $C_k$ by $u_k$. 

The {\it Simis cone} of $I=I(\mathcal{C})$ is the
rational polyhedral cone: 
$$
{\rm Cn}(I)=H_{e_1}^+\cap\cdots\cap H_{e_{n+1}}^+\cap
H_{(u_1,-1)}^+\cap\cdots\cap H_{(u_s,-1)}^+.
$$
Here $H_{a}^+$ denotes the closed halfspace 
$H_a^+=\{x\vert\, \langle
x,a\rangle\geq 0\}$ and $H_a$ stands for the hyperplane through the origin with normal
vector $a$. Simis cones were 
introduced in \cite{normali} to study symbolic Rees algebras of
square-free monomial ideals.  The term {\it Simis cone} is intended to do homage to Aron
Simis \cite{aron-hoyos,simis-ulrich,ITG}. The Simis cone is a pointed rational polyhedral
cone. By \cite[Theorem~16.4]{Schr} there is a unique minimal finite set of
integral vectors 
$$
\mathcal{H}=\{h_1,\ldots,h_r\}\subset\mathbb{Z}^{n+1}
$$
such that 
$\mathbb{Z}^{n+1}\cap \mathbb{R}_+\mathcal{H}=\mathbb{N}\mathcal{H}$ and
${\rm Cn}(I)=\mathbb{R}_+\mathcal{H}$ (minimal relative to
taking subsets), where  $\mathbb{R}_+\mathcal{H}$ denotes the
cone generated by $\mathcal{H}$ consisting of all linear combinations of
$\mathcal{H}$ with non-negative real coefficients and
$\mathbb{N}\mathcal{H}$ denotes the semigroup generated by
$\mathcal{H}$ consisting of all linear combinations of
$\mathcal{H}$ with coefficients in $\mathbb{N}$. The set $\mathcal{H}$
is called the {\it Hilbert basis\/} of ${\rm Cn}(I)$. 
The Hilbert basis of
${\rm Cn}(I)$ has the following useful description.

\begin{theorem}{\rm \cite[p. 233]{Schr}}\label{hb-description}
$\mathcal{H}$ is the set of all integral vectors $0\neq h\in {\rm Cn}(I)$ such
that $h$ is not the sum of two other non-zero integral vectors in
${\rm Cn}(I)$. 
\end{theorem}

\begin{corollary}\label{april14-09} Let $\mathcal{H}$ be the Hilbert basis of 
${\rm Cn}(I)$. Then
\begin{eqnarray}
\mathcal{H}&=&\{(a,b)\vert\, x^at^b\mbox{ is a
minimal generator of }R_s(I)\}\\
&=&\{(a,\alpha_0(\mathcal{C}))\vert\, \mathcal{C}^a\mbox{ is an
indecomposable parallelization of }\mathcal{C}\}
\end{eqnarray}
and $R_s(I)$ is equal to the semigroup ring $K[\mathbb{N}{\mathcal H}]$ of 
$\mathbb{N}{\mathcal H}$. 
\end{corollary}

\begin{proof} The first equality follows from Lemma~\ref{april8-09}
and Theorem~\ref{hb-description}.
The second equality follows from Theorem~\ref{irr-b-cover-char}. The
equality $K[\mathbb{N}{\mathcal H}]=\mathbb{N}{\mathcal H}$ was first
observed in \cite[Theorem~3.5]{normali}.
\end{proof}
This result is interesting because it 
allows to compute all indecomposable parallelizations of $\mathcal{C}$ and
all indecomposable induced subclutters of $\mathcal{C}$ using Hilbert
bases. In particular, as is seen in
Corollary~\ref{method-irred-sub}, we can use this result to decide whether any
given graph or clutter is indecomposable (see Example~\ref{example-reyes}).

The indecomposable subclutters can be computed using the next 
consequence of Corollary~\ref{april14-09}.

\begin{corollary}\label{method-irred-sub} Let $\mathcal{C}$ be a
clutter and let
$\alpha=(a_1,\ldots,a_n,b)$ be a vector in $\{0,1\}^n\times\mathbb{N}$. Then 
$\alpha$ is in the Hilbert basis of ${\rm
Cn}(I(\mathcal{C}))$ if and only if the induced subclutter 
$H=\mathcal{C}[\{x_i\vert\, a_i=1\}]$ is indecomposable with
$b=\alpha_0(H)$.
\end{corollary}

\begin{example}\label{example-reyes} Consider the
graph $G$ 
shown below.
Let $I$ be the edge ideal of $G$ and let $\mathcal{H}$ be the Hilbert
basis of ${\rm Cn}(I)$. Using Corollary~\ref{april14-09}, together with 
{\sc Normaliz} \cite{normaliz2}, it is seen that $G$ has exactly 61
indecomposable parallelizations and 49 indecomposable subgraphs. Since 
$\alpha_0(G)=6$ and the vector $(1,\ldots,1,6)$ is not in
$\mathcal{H}$ we obtain that $G$ is a decomposable graph. 

\setlength{\unitlength}{.035cm}
\thicklines
\begin{picture}(10,30)(-100,20)
\put(-30,0){\circle*{4.2}}
\put(30,0){\circle*{4.2}}
\put(0,20){\circle*{4.2}}
\put(-20,-20){\circle*{4.2}}
\put(20,-20){\circle*{4.2}}
\put(-30,0){\line(3,2){30}}
\put(30,0){\line(-3,2){30}}
\put(-30,0){\line(1,-2){10}}
\put(-20,-20){\line(1,0){40}}
\put(20,-20){\line(1,2){10}}
\put(70,0){\circle*{4.2}}
\put(130,0){\circle*{4.2}}
\put(100,20){\circle*{4.2}}
\put(80,-20){\circle*{4.2}}
\put(120,-20){\circle*{4.2}}
\put(70,0){\line(3,2){30}}
\put(130,0){\line(-3,2){30}}
\put(70,0){\line(1,-2){10}}
\put(80,-20){\line(1,0){40}}
\put(120,-20){\line(1,2){10}}
\put(-43,0){$x_5$}
\put(35,-5){$x_6$}
\put(0,25){$x_9$}
\put(-33,-17){$x_1$}
\put(5,-15){$x_2$}
\put(75,-2){$x_7$}
\put(135,0){$x_8$}
\put(100,25){$x_{10}$}
\put(80,-15){$x_3$}
\put(125,-20){$x_4$}
\qbezier(-20,-20)(-10,-40)(80,-20)
\qbezier(-20,-20)(-10,-50)(120,-20)
\qbezier(20,-20)(120,-20)(80,-20)
\qbezier(20,-20)(80,10)(120,-20)
\qbezier(0,20)(60,0)(70,0)
\qbezier(0,20)(90,0)(130,0)
\qbezier(100,20)(20,0)(30,0)
\qbezier(100,20)(-20,0)(-30,0)
\put(-5,-50){Fig. 4. Decomposable graph $G$ }
\end{picture}

\vspace{2.8cm}

The vector 
$a=(1,\ldots,1,2,7)$ is in $\mathcal{H}$, i.e., $G^{(1,\ldots,1,2)}$
is indecomposable and has covering number $7$. 

\setlength{\unitlength}{.035cm}
\thicklines
\begin{picture}(10,30)(-100,20)
\put(-30,0){\circle*{4.2}}
\put(30,0){\circle*{4.2}}
\put(0,20){\circle*{4.2}}
\put(-20,-20){\circle*{4.2}}
\put(20,-20){\circle*{4.2}}
\put(-30,0){\line(3,2){30}}
\put(30,0){\line(-3,2){30}}
\put(-30,0){\line(1,-2){10}}
\put(-20,-20){\line(1,0){40}}
\put(20,-20){\line(1,2){10}}
\put(50,30){\circle*{4.2}}
\put(70,0){\circle*{4.2}}
\put(130,0){\circle*{4.2}}
\put(100,20){\circle*{4.2}}
\put(80,-20){\circle*{4.2}}
\put(120,-20){\circle*{4.2}}
\put(70,0){\line(3,2){30}}
\put(130,0){\line(-3,2){30}}
\put(70,0){\line(1,-2){10}}
\put(80,-20){\line(1,0){40}}
\put(120,-20){\line(1,2){10}}
\qbezier(-20,-20)(-10,-40)(80,-20)
\qbezier(-20,-20)(-10,-50)(120,-20)
\qbezier(20,-20)(120,-20)(80,-20)
\qbezier(20,-20)(80,10)(120,-20)
\qbezier(0,20)(60,0)(70,0)
\qbezier(0,20)(90,0)(130,0)
\qbezier(100,20)(20,0)(30,0)
\qbezier(100,20)(-20,0)(-30,0)
\qbezier(50,30)(0,10)(-30,0)
\put(50,30){\line(-2,-3){20}}
\put(50,30){\line(2,-3){20}}
\qbezier(50,30)(120,0)(130,0)
\put(-43,0){$x_5$}
\put(35,-5){$x_6$}
\put(0,25){$x_9$}
\put(-33,-17){$x_1$}
\put(5,-15){$x_2$}
\put(75,-2){$x_7$}
\put(135,0){$x_8$}
\put(100,25){$x_{10}$}
\put(80,-15){$x_3$}
\put(125,-20){$x_4$}
\put(55,32){$x_{10}'$}
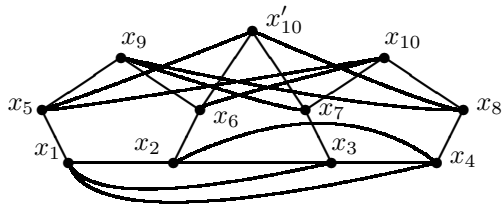
\put(-15,-50){Fig. 5. Indecomposable graph $G^{(1,\ldots,1,2)}$}
\end{picture}
\end{example}
\vspace{2.5cm}

The next result, together with Corollary~\ref{method-irred-sub}, 
allows to locate all {\it induced odd cycles\/} ({\it odd
holes\/}) and all {\it induced 
complements of odd cycles\/} ({\it odd antiholes}). 

\begin{lemma}\label{antiholes-are-irred} Let $C_n=\{x_1,\ldots,x_n\}$
be a cycle. $(\mathrm{a})$ If $n\geq 5$ is odd, then the complement
$C_n'$ of $C_n$ is an indecomposable graph, $(\mathrm{b})$ 
if $n$ is odd, then $C_n$ is an indecomposable cycle, and $(\mathrm{c})$
any complete graph is indecomposable.
\end{lemma}

\begin{proof} (a) Assume that $G=C_n'$ is decomposable. Then there are disjoint sets
$X_1,X_2$ such that $V(G)=X_1\cup X_2$ and
$\alpha_0(G)=\alpha_0(G[X_1])+\alpha_0(G[X_2])$. Since $\beta_0(G)=2$,
it is seen that $G[X_i]$ is a complete graph for $i=1,2$. We may
assume that $x_1\in X_1$. Then $x_2$ must be in $X_2$, otherwise
$\{x_1,x_2\}$ is an edge of $G[X_1]$, a contradiction. By induction 
it follows that $x_1,x_3,x_5,\ldots,x_n$ are in $X_1$. Consequently 
$\{x_1,x_n\}$ is an edge of $G[X_1]$, a contradiction. Thus $G$ is
indecomposable. (b) This was observed in \cite{harary-plummer}. (c)
Follows readily from the fact that the covering 
number of a complete graph in $r$ vertices is $r-1$. 
\end{proof}

\begin{example} Consider the graph $G$ of Fig. 6, where vertices 
are labeled with $i$ instead of $x_i$. Using
Corollary~\ref{april14-09}, together with {\sc Normaliz} \cite{normaliz2}, it is
seen that $G$ has exactly $21$ indecomposable parallelizations, 
$20$ of which correspond to indecomposable subgraphs. Apart from the seven
vertices, 
the nine edges, one triangle and three pentagons, the only
indecomposable parallelization of $G$ which is not a subgraph is the
duplication shown in Fig 7.
$$
\begin{array}{cccc}
\setlength{\unitlength}{.04cm} \thicklines
\begin{picture}(0,-50)(10,20)
\put(0,0){\circle*{4.2}} \put(60,0){\circle*{4.2}} \put(0,30){\circle*{4.2}}
\put(30,60){\circle*{4.2}} \put(60,30){\circle*{4.2}}\put(30,30){\circle*{4.2}}
\put(30,15){\circle*{4.2}}

\put(0,0){\line(1,0){60}}\put(0,0){\line(0,1){30}}\put(0,0){\line(2,1){30}}
\put(60,0){\line(0,1){30}}
\put(0,30){\line(1,1){30}}\put(60,30){\line(-1,1){30}}\put(30,15){\line(0,1){15}}
\put(30,30){\line(0,1){30}}\put(60,0){\line(-2,1){30}}

\newcommand{\lb}[1]{\tiny $#1$}
\put(-6,0){\lb{4}} \put(64,0){\lb{3}} \put(-6,28){\lb{5}} \put(29,63){\lb{1}}
\put(64,28){\lb{2}}\put(24,28){\lb{6}} \put(24,15){\lb{7}}
\put(-26,-20){Fig. 6. Decomposable graph $G$}
\end{picture}
& \ \ 
&
\setlength{\unitlength}{.04cm} \thicklines
\begin{picture}(10,50)(-150,20)
\put(0,0){\circle*{4.2}} 
\put(60,0){\circle*{4.2}} \put(0,30){\circle*{4.2}}
\put(30,60){\circle*{4.2}} \put(60,30){\circle*{4.2}}\put(30,30){\circle*{4.2}}
\put(30,15){\circle*{4.2}}\put(20,40){\circle*{4.2}}

\put(0,0){\line(1,0){60}}\put(0,0){\line(0,1){30}}\put(0,0){\line(2,1){30}}
\put(60,0){\line(0,1){30}}
\put(0,30){\line(1,1){30}}\put(60,30){\line(-1,1){30}}\put(30,15){\line(0,1){15}}
\put(30,30){\line(0,1){30}}\put(60,0){\line(-2,1){30}}\put(20,40){\line(-2,-1){20}}
\put(20,40){\line(1,-1){10}}\put(20,40){\line(4,-1){40}}

\newcommand{\lb}[1]{\tiny $#1$}
\put(-6,0){\lb{4}} \put(64,0){\lb{3}} \put(-6,28){\lb{5}} \put(29,63){\lb{1}}
\put(64,28){\lb{2}}\put(24,28){\lb{6}} \put(24,15){\lb{7}}\put(22,42){\lb{1'}}
\put(-68,-20){Fig. 7. Indecomposable graph $G^{(2,1,1,1,1,1,1)}$}
\end{picture}
&\ \ \ \ \ \ \ \ \ \ \ \ \ \ \ \ \ \ \ \ \ \ \ \ \ \ \ \ \ \ \ \ \ \ 
\ \ \ \ \ \ \ \ \ \ \ \ \ \ \ \  \ \ \ \ \ \ \ \ \ \ \ \ \ \ \ \ 
\end{array}
$$
\end{example}

\vspace{1.5cm}

\begin{example} Consider the graph 
$G$ of Fig. 8. Using Corollary~\ref{april14-09} and {\sc Normaliz}
\cite{normaliz2}, it is seen that $G$ has exactly $103$ indecomposable
parallelizations, $92$ of which correspond to indecomposable subgraphs.
The only indecomposable parallelization $G^a$ which do not delete vertices is that
obtained by duplication of the five outer vertices, i.e., 
$a=(2,  2,  2, 2,  2,  1, 1,  1, 1 , 1)$ and
$\alpha_0(G^a)=11$.

\setlength{\unitlength}{.035cm} \thicklines
\begin{picture}(10,95)(-100,20)
 \put(0,0){\circle*{4.2}} \put(60,0){\circle*{4.2}} \put(30,100){\circle*{4.2}}
 \put(20,20){\circle*{4.2}} \put(40,20){\circle*{4.2}}\put(10,50){\circle*{4.2}}
 \put(50,50){\circle*{4.2}}\put(30,70){\circle*{4.2}}\put(-20,50){\circle*{4.2}}
 \put(80,50){\circle*{4.2}}

 \put(0,0){\line(1,0){60}} \put(0,0){\line(-2,5){20}} \put(60,0){\line(2,5){20}}
 \put(-20,50){\line(1,1){50}}
 \put(80,50){\line(-1,1){50}}
 \put(20,20){\line(1,1){30}}
 \put(20,20){\line(1,5){10}}
 \put(10,50){\line(1,0){40}}
 \put(10,50){\line(1,-1){30}}
 \put(40,20){\line(-1,5){10}}

 \put(0,0){\line(1,1){20}}
 \put(60,0){\line(-1,1){20}}
 \put(-20,50){\line(1,0){30}}
 \put(50,50){\line(1,0){30}}
 \put(30,70){\line(0,1){30}}

 \put(0,0){\line(1,5){10}}
 \put(0,0){\line(2,1){40}}
 \put(20,20){\line(-4,3){40}}
 \put(60,0){\line(-2,1){40}}
 \put(60,0){\line(-1,5){10}}
 \put(40,20){\line(4,3){40}}
 \put(-20,50){\line(5,2){50}}
 \put(10,50){\line(2,5){20}}
 \put(50,50){\line(-2,5){20}}
 \put(80,50){\line(-5,2){50}}

 \newcommand{\lb}[1]{\tiny $#1$}
 \put(-8,0){\lb{4}} \put(66,0){\lb{3}} \put(30,105){\lb{1}}
 \put(16,24){\lb{9}} \put(41,24){\lb{8}}\put(2,52){\lb{10}}
 \put(52,52){\lb{7}}\put(24,72){\lb{6}}\put(-26,50){\lb{5}}
 \put(84,50){\lb{2}}
\put(-26,-20){Fig. 8. Decomposable graph $G$}
\end{picture}
\end{example}

\vspace{1.5cm} 

\section{Symbolic Rees algebras and perfect
graphs}\label{perfect-graphs}

We now turn our attention to the indecomposability of graphs and its
connection with the theory of perfect graphs. 
Examples of indecomposable graphs include complete graphs, odd cycles, and
complements of odd cycles of length at least $5$ (see Lemma~\ref{antiholes-are-irred}). 

Let us recall the notion of a perfect graph that was introduced 
by Berge \cite[Chapter 16]{berge-graphs-hypergraphs}. 
A {\it colouring\/} of the vertices of a graph $G$ is an assignment
of colours to the vertices of $G$ in such a way that adjacent vertices
have distinct colours. The {\it chromatic number\/} of $G$ is the
minimal number of colours in a colouring of $G$. 
A graph is {\it perfect\/} if for every induced subgraph $H$, the
chromatic 
number of $H$ equals the size of the largest complete subgraph 
of $H$. We refer to \cite {cornu-book,golumbic,Schr2} for the theory
of perfect graphs.  

The next result shows that indecomposable
graphs occur naturally in the theory of perfect graphs. 

\begin{proposition}{\rm \cite[Proposition~2.13]{covers}}
\label{perfect-char-parallel} 
A graph $G$ is perfect if and only if the indecomposable parallelizations
of $G$ are exactly the complete subgraphs of $G$ 
\end{proposition}

Let $G$ be a graph. We denote a complete
subgraph of $G$ with $r$ vertices  
by ${\mathcal K}_r$. The empty
set is regarded as an independent set of vertices whose characteristic vector is the
zero vector. A {\it clique\/} of $G$ is a subset of the set
of vertices that induces 
a complete subgraph. The {\it support\/} of a monomial 
$x^a=x_1^{a_1}\cdots x_n^{a_n}$, denoted by ${\rm supp}(x^a)$, is the
set ${\rm supp}(x^a)= \{x_i\, |\, 
a_i>0\}$. If $a_i\in\{0,1\}$ for all $i$, $x^a$ is called a {\it
square-free\/} monomial.  

The next major result shows that the symbolic Rees algebra of the edge
ideal of a perfect 
graph $G$ is completely determined by the cliques of $G$. This
was first shown in \cite{perfect} using polyhedral geometry.

\begin{corollary}{\cite[Corollary~3.3]{perfect}}\label{coro2-irr-b-cover-char}
If $G$ is a perfect graph, then 
$$
R_s(I(G))=K[\{x^at^b\vert\, 
x^a\mbox{ is square-free };\, 
G[{\rm supp}(x^a)]={\mathcal K}_{b+1}\}].
$$
\end{corollary}

\begin{proof} Let $x^at^b$ be a minimal of $R_s(I(G))$. By
Theorem~\ref{irr-b-cover-char} $G^a$ is an 
indecomposable graph and $b=\alpha_0(G^a)$. As $G$ is perfect, by
Proposition~\ref{perfect-char-parallel}, we obtain that $G^a$ is a
complete subgraph of $G$ with $b+1$ vertices. 
\end{proof} 

Since complete graphs are perfect, an immediate consequence is: 

\begin{corollary}{\rm \cite{bahiano}}
If $G$ is a complete graph, then 
$$
R_s(I(G))=K[\{x^at^b\vert\, 
x^a\mbox{ is square-free };\, \deg(x^a)={b+1}\}].
$$
\end{corollary}

\bibliographystyle{amsplain}

\end{document}